\documentclass[11pt]{article}

\usepackage{amssymb}
\usepackage{amsmath}

\newtheorem{theorem}{Theorem}[section]
\newtheorem{lemma}[theorem]{Lemma}

\newtheorem{corollary}[theorem]{Corollary}

\numberwithin{equation}{section}

\newenvironment{proof}{\vspace{1.5mm}\newline \bf Proof. \rm}{$\Box$}

\begin{document}

\title{Using the smoothness of $p-1$ for computing roots modulo $p$}

\author{Bartosz \'Zra{\l}ek
\footnote{Institute of Mathematics, Polish Academy of Sciences, 00-956 Warsaw, Poland.
email:~\textit{b.zralek@impan.gov.pl}}}

\maketitle





\begin{abstract}
We prove, without recourse to the Extended Riemann Hypothesis, that the projection modulo $p$ 
of any \textit{prefixed} polynomial with integer coefficients can be completely factored in deterministic polynomial time if  $p-1$ has a $(\ln p)^{O(1)}$-smooth divisor exceeding $(p-1)^{\frac{1}{2}+\delta}$ 
for some arbitrary small $\delta$. We also address the issue of computing roots modulo $p$ in
deterministic time.
\end{abstract}

\section{Introduction}

Factoring polynomials over finite fields in deterministic polynomial time is a long-standing open problem of computational number theory. The most important results obtained, though partial, are now classic. Berlekamp \cite{berlekamp} was the first to devise a general deterministic algorithm for this problem;
its running time bound $p(d\ln p)^{O(1)}$, where $p$ is the characteristic of the finite field and $d$ the degree of the polynomial to be factored, can be seen as polynomial only if $p$ is fixed. A better and so far best time bound $p^\frac{1}{2}(d\ln p)^{O(1)}$ is achieved by an algorithm of Shoup \cite{shoup1}. There are also algorithms with running time bound of the form $(d\ln p)^{O(1)}$, such as the Cantor-Zassenhaus algorithm \cite{cz} (actually dating back to Legendre), but these use randomness in an essential way. In this article we pursue an approach developed by von zur Gathen
\cite{gathen} and R\'onyai \cite{ronyai}, which consists in taking advantage of the multiplicative
structure of $p-1$. Shoup \cite{shoup2} refined the algorithms of von zur Gathen and R\'onyai, improving the running time bound from $P^+(p-1)(d\ln p)^{O(1)}$ to $P^+(p-1)^\frac{1}{2}(d\ln p)^{O(1)}$, where $P^+(p-1)$ is the largest prime factor of $p-1$. These three algorithms are deterministic, however in their proofs of correctness the Extended Riemann Hypothesis must be assumed. We prove, without recourse to
the ERH, that the projection modulo $p$ of any \textit{prefixed} polynomial with integer coefficients can
be completely factored in deterministic polynomial time if  $p-1$ has a $(\ln p)^{O(1)}$-smooth
divisor exceeding $(p-1)^{\frac{1}{2}+\delta}$ for some arbitrary small $\delta$.  

\begin{theorem}
\label{main}
Let $f$ be an irreducible polynomial of degree $d$ in $\mathbb{Z}[X]$, with leading coefficient $l$,
and $\tilde{f}$ the polynomial of $\mathbb{Z}[Y]$ defined by $\tilde{f}(Y)=l^{d-1}f(\frac{Y}{l})$. 
Let $h$ be the class number of $\mathbb{Q}(\theta)$, where $\theta$ is any complex root of $\tilde{f}$. 
Let $p$ be a prime and $q$ the least prime such that the $q$-smooth part
$S$ of $p-1$ is no less than $(p-1)^{\frac{1}{2}+\delta}$ for some $\delta>0$.
Then the complete factorization of $f$ modulo $p$ can be found in 
$O_{\theta, c, \delta, \epsilon}((q^\frac{1}{2}\ln q+\ln p)\ln^{2+\gamma+ch} p)$ deterministic time, where
$\gamma$ is any positive number and $c, \epsilon$ positive numbers satisfying 
$\frac{d}{c}+\epsilon<2\delta$.
\end{theorem}

Let us emphasize that the above big-$O$ constant depends severely on $\theta$ as to include the time
of some necessary precomputations  in the number field $\mathbb{Q}(\theta)$. For a thorough 
exposition of constructive algebraic number theory we refer the reader to \cite{pohst}.
In a sense, the idea of fixing the polynomial $f$ rather than the prime $p$ is at the opposite of
Berlekamp's algorithm. Actually, it has been already considered by Schoof \cite{schoof} in the case when $f=X^2-a$ for $a$ an integer, and by Pila \cite{pila} when $f=X^s-1$ for $s$ a prime dividing 
$p-1$. Both authors gave unconditional algorithms that factor the corresponding
polynomial $f$ modulo $p$ in respectively $O_a(\ln^{6+\epsilon} p)$ and $(\ln p)^{O_s(1)}$
deterministic time. We address more generally the issue of computing roots modulo $p$ in the ensuing corollary. Our proofs are mainly based upon algebraic number theory, whereas Schoof's and Pila's rely on heavy machinery of algebraic geometry and this is of independent interest.

\begin{corollary}
\label{pierwiastki}
Let $p$ be a prime and $q$ the least prime such that the $q$-smooth part of $p-1$ is no less
than $(p-1)^{\frac{1}{2}+\delta}$ for some $\delta>0$. Let $n$ be a positive integer. Suppose that the integer $a$ is an $n$-th power residue modulo $p$. Then all the $n$-th roots of $a$ modulo $p$ can be computed in $O_{a, n, c, \delta, \epsilon}((q^\frac{1}{2}\ln q+\ln p)\ln^{2+\gamma+cH} p)$ deterministic time, where $\gamma$ is any positive number, $c, \epsilon$ positive numbers satisfying 
$\frac{n}{c}+\epsilon<2\delta$, and $H$ the largest integer among the class numbers of the fields 
$\mathbb{Q}(\theta)$, $\theta$ running through the complex $n$-th roots of $a$.     
\end{corollary}

In the cases covered by Schoof or Pila and such that the corresponding integer $H$ is equal to one
the stated running time bound is slightly better for a sparse, but infinite set of primes $p$.
It is worth noting that our technique combined with the result of Pila and an observation of Tsz-Wo Sze
\cite{sze} leads for $n$ an odd prime to a stronger theorem than corollary \ref{pierwiastki} 
(see last remark of section \ref{remarks}).


\section{Notation}

In all that follows, $f$ is a \textit{fixed}, irreducible polynomial of degree $d$
in $\mathbb{Z}[X]$, with leading coefficient $l$, and $\tilde{f}$ the corresponding monic, irreducible
polynomial $l^{d-1}f(\frac{Y}{l})$ of degree $d$ in $\mathbb{Z}[Y]$.
The number field $K$ is the extension of $\mathbb{Q}$ by a complex root $\theta$ of the
polynomial $\tilde{f}$. The class number of $K$ is $h$, its ring of integers - $\mathcal{O}_K$.
A \textit{fixed}, integral basis $\omega=(\omega_1,\ldots,\omega_d)$ of $\mathcal{O}_K$,
as well as a \textit{fixed}, finite set $\mathcal{U}$ of generators of the group of units $\mathcal{O}_K^*$
are given. The ideals of $\mathcal{O}_K$ that we consider are always supposed to be nonzero.
The norm $N(I)$ of an ideal in $\mathcal{O}_K$ is the cardinality of $\mathcal{O}_K/I$.
We let $\psi_K(x, y)$, respectively $\tilde{\psi}_K(x, y)$, be the number of ideals, respectively
principal ideals, of $\mathcal{O}_K$ with norm at most $x$ that can be written as a product of prime ideals, respectively principal ideals, of $\mathcal{O}_K$ with norm at most $y$.
\\
The letter $p$ denotes an odd prime number. 
For $g\in\mathbb{Z}_p[Y]$, by $R_g$ we mean the quotient ring $\mathbb{Z}_p[Y]/(g)$
and by $R_g^*$ its multiplicative group. If the commutative group $G$ is a direct sum of
two subgroups $G_1$, $G_2$, and $\mathcal{G}$ is a subset of $G$ then the symbol 
$\langle\mathcal{G}\rangle_G$, respectively $\langle\mathcal{G}\rangle_{G_1}$, stands for the subgroup of $G$ generated by $\mathcal{G}$, respectively the subgroup of $G/G_2$ generated 
by the cosets $gG_2$, $g\in\mathcal{G}$.

\section{Auxiliary results}

We will seek to compute the factorization of $\tilde{f}$ modulo $p$; it gives the factorization of $f$ by a change of variable whenever $p$ does not divide $l$. In general, we can assume that the prime $p$ exceeds any given constant, since for small $p$ factoring in $\mathbb{Z}_p[X]$ is "easy".

The algorithm of Fellows and Koblitz \cite{fellows} for proving the primality of an integer 
$n$ or also the deterministic version of Pollard's $p-1$ method \cite{zralek} for factoring $n$, perform certain operations on a "small" subset $\mathcal{B}$ of $\mathbb{Z}$ that generates modulo $n$ a "large" multiplicative semigroup $\mathcal{S}$. In fact, $\mathcal{B}$ can be chosen 
there as the set of prime numbers not exceeding $\ln^2 n$; the semigroup $\mathcal{S}$ has then at least $\psi(n, \ln^2 n)$ elements, where $\psi$ is the de Bruijn function. Here similarly, to factor $\tilde{f}$ modulo $p$ we will construct a small subset of $\mathbb{Z}_p[Y]/(\tilde{f})$ generating a large multiplicative semigroup. By the following lemma, the latter task amounts to exhibiting a suitable subset of $\mathcal{O}_K$, at least if $p$ is sufficiently large.

\begin{lemma}
\label{izomorfizm}
Assume that $p$ does not divide the index $[\mathcal{O}_K:\mathbb{Z}[\theta]]$. Then 
$\theta\mapsto Y$ induces an isomorphism 
$\kappa:~\mathcal{O}_K/(p)\rightarrow\mathbb{Z}_p[Y]/(\tilde{f})$.
\end{lemma}

Unlike $\mathbb{Z}$ however, $\mathcal{O}_K$ is not a unique
factorization domain (unless $h=1$). It is still a Dedekind domain and just as $\psi$ measures the smoothness of integers in $\mathbb{Z}$, the function $\psi_K$ measures the smoothness of ideals
in $\mathcal{O}_K$. The next theorem generalizes in this sense a result of Canfield et al. \cite{canfield}.  

\begin{theorem}[Moree, Stewart]
\label{moree}
There is an effective, positive constant $c_1=c_1(K)$ such that for $x, y\ge 1$
and $u:=\frac{\ln x}{\ln y}\ge 3$ we have
\begin{equation*}
\psi_K(x, y)\ge x\exp
\left[-u\left\{\ln(u\ln u)-1+\frac{\ln\ln u-1}{\ln u}+c_1\left(\frac{\ln\ln u}{\ln u}\right)^2\right\}\right].
\end{equation*}
\end{theorem}

The generator of a principal ideal in $\mathcal{O}_K$ is defined up to a multiplicative unit of
$\mathcal{O}_K$, so working with principal ideals, rather than general ones, is pretty much like
working with integers. That is why we give via $\psi_K$ a lower bound for the function $\tilde{\psi}_K$ counting the number of "smooth" principal ideals.
 
\begin{lemma}
\label{ideal}
There is an effective, positive constant $c_2=c_2(K)$ such that 
$\tilde{\psi}_K(x, y)\ge~${\Large$\frac{1}{h}$}$\psi_K(c_2x, y^\frac{1}{h})$ for $y\ge c_2^{-h}$.
\begin{proof}
Let $I_1,\ldots,I_h$ be a set of representatives for the class group of $K$ whose norms are bounded
above by the Minkowski bound $M_K$. We will prove that the lemma holds with $c_2=M_K^{-1}$.
Define $\psi_K'(x, y)$ as the number of principal ideals of $\mathcal{O}_K$ with norm at most $x$
that split as a product of prime ideals of $\mathcal{O}_K$ with norm at most $y$.
\\
Let $J$ be an ideal counted by $\psi_K(M_K^{-1}x, y^\frac{1}{h})$. There exists a $k$, $1\le k\le h$,
such that $JI_k$ is principal. Suppose that $y^\frac{1}{h}\ge M_K$, i.e. $y\ge M_K^h$. Then $JI_k$
is counted by $\psi_K'(x, y^\frac{1}{h})$. Moreover, any ideal counted by $\psi_K'(x, y^\frac{1}{h})$
can be written in at most $h$ ways as $JI_k$, where $J$ is counted by 
$\psi_K(M_K^{-1}x, y^\frac{1}{h})$ and $1\le k\le h$. Consequently,
{\Large$\frac{1}{h}$}$\psi_K(M_K^{-1}x, y^\frac{1}{h})\le\psi_K'(x, y^\frac{1}{h})$.
\\
Assume that the principal ideal $I$ of $\mathcal{O}_K$ is a product of $m$ prime ideals of
$\mathcal{O}_K$ with norm at most $y^\frac{1}{h}$. It is easy to show by induction on $m$ that
$I$ is a product of principal ideals of $\mathcal{O}_K$ with norm at most $y$. Just use the fact
that every product of at least $h$ ideals of $\mathcal{O}_K$ contains a principal factor.
Therefore any ideal counted by  $\psi_K'(x, y^\frac{1}{h})$ is also counted by $\tilde{\psi}_K(x, y)$,
hence $\psi_K'(x, y^\frac{1}{h})\le\tilde{\psi}_K(x, y)$.
\end{proof}
\end{lemma}

It becomes apparent that if we let $\mathcal{B}$ be the union of $\mathcal{U}$ and a set $\mathcal{A}$
containing pairwise non-associate integers with small norm then it should generate modulo $p$ a relatively large multiplicative semigroup $\mathcal{S}$. Nevertheless, three problems arise. Is 
$\mathcal{S}$ indeed large though reduction modulo $p$? Can the suitable set $\mathcal{A}$ be small and easy to find? The ensuing theorem helps to answer these questions positively.   

\begin{theorem}[Fincke, Pohst]
\label{pohst}
There is an effective, positive constant $c_3=c_3(K, \omega)$ such that for any 
$\eta\in\mathcal{O}_K\setminus\{0\}$
there exists $\tilde{\eta}\in\mathcal{O}_K$ generating the same ideal as $\eta$ and whose coordinates $a_i$ in the basis $\omega$ satisfy $|a_i|\le c_3N((\eta))^\frac{1}{d}$.
\begin{proof}
Combine the equations (3.5b), chapter 5, and (4.3f), chapter 6 of \cite{pohst}.
\end{proof}
\end{theorem}

We now summarize rigorously the above informal discussion. Actually, we show more: for any
$g$ dividing $f$ modulo $p$, a set $\mathcal{B}_g$ derived from $\mathcal{B}$ generates a large multiplicative semigroup in $R_g$.
  
\begin{lemma}
\label{baza}
Suppose that the polynomial $g$ of degree $d'$ divides $\tilde{f}$ modulo $p$. Let $p$ and $\kappa$ 
be as in lemma \ref{izomorfizm}. Let $\pi$ and $\pi_g$ be the projections 
$\mathcal{O}_K\rightarrow\mathcal{O}_K/(p)$ and 
$\mathbb{Z}_p[Y]/(\tilde{f})\rightarrow\mathbb{Z}_p[Y]/(g)$ respectively.
Fix $c>0$ and define
$\mathcal{A}=\{a_1\omega_1+\ldots+a_d\omega_d:~
a_i\in\mathbb{Z},~|a_i|\le c_3\ln^{\frac{ch}{d}}p,~1\le i\le d\}$,
$\mathcal{S}=\{v\cdot\alpha_1\cdot\ldots\cdot\alpha_m:~v\in\mathcal{O}_K^*,~m\in\mathbb{N},~
\alpha_i\in\mathcal{A},~1\le i\le m\}$.
Then $\#\pi_g\kappa\pi(\mathcal{S})> p^{d'-\frac{d}{c}-\epsilon}$ for any $\epsilon>0$ and 
$p\ge p_0$, $p_0=p_0(c, c_1, c_2, c_3, \epsilon)$.
\begin{proof}
Let $\mathcal{T}=\mathcal{S}\cap
\{a_1\omega_1+\ldots+a_d\omega_d:~a_i\in\mathbb{Z},~|a_i|\le\frac{p}{2},~1\le i\le d\}$. It is
sufficient to prove that the desired inequality holds with $\mathcal{S}$ replaced by $\mathcal{T}$.
We invoke theorem \ref{pohst} to get
$\#\mathcal{T}\ge\tilde{\psi}_K((\frac{p}{2c_3})^d, \ln^{ch}p)$.
By lemma \ref{ideal}, if $p$ is large enough the latter expression is no less than
{\Large$\frac{1}{h}$}$\psi_K(c_2(\frac{p}{2c_3})^d, \ln^c p)$. This in turn is greater than
$p^{d-\frac{d}{c}-\epsilon}$ for any $\epsilon>0$ and sufficiently large $p$, by theorem \ref{moree}.
Thus $\#\mathcal{T}>p^{d-\frac{d}{c}-\epsilon}$ if $p$ exceeds some constant $p_0$ depending
upon $c, c_1, c_2, c_3$ and $\epsilon$. Assume that it does. As $p>2$, we have 
$\#\pi(\mathcal{T})=\#\mathcal{T}$. Furthermore, $\kappa$ is an isomorphism, hence
$\#\kappa\pi(\mathcal{T})=\#\pi(\mathcal{T})$. Finally, $\pi_g$ is a surjective homomorphism,
so the preimage under $\pi_g$ of any element of $\mathbb{Z}_p[Y]/(g)$ has
$\#\ker\pi_g=p^{d-d'}$ elements. It follows that 
$p^{d-d'}\cdot\#\pi_g\kappa\pi(\mathcal{T})\ge\#\kappa\pi(\mathcal{T})=\#\mathcal{T}>
p^{d-\frac{d}{c}-\epsilon}$. Therefore $\#\pi_g\kappa\pi(\mathcal{T})>p^{d'-\frac{d}{c}-\epsilon}$.
\end{proof}
\end{lemma}

Assume that $g$ is a product of at least two distinct, degree $e$ irreducible factors of $f$ modulo $p$.
Either the set $\mathcal{B}_g$ mentioned above is not contained in $R_g^*\cup\{0\}$, or 
$\mathcal{G}:=\{b^\frac{p^e-1}{p-1}:~b\in\mathcal{B}_g\setminus\{0\}\}$ generates a large subgroup
of  $\{a\in R_g^*:~a^{p-1}=1\}$ and thus should not be cyclic. The latter case is dealt with an extension
of the Pohlig-Hellman algorithm \cite{pohlig} for computing discrete logarithms.

\begin{theorem}
\label{PH}
Let $g$ be a polynomial of degree $d'$ in $\mathbb{Z}_p[Y]$ and
$G$ the group $\{a\in R_g^*:~a^{p-1}=1\}$. 
Write $G=G_1\oplus G_2$, $(\#G_1,\#G_2)=1$.
Suppose that we are given $g$, $\#G_1$ and a subset $\mathcal{G}$ of $G$ such that 
$\langle\mathcal{G}\rangle_{G_1}$ is not cyclic. Then we can find a nontrivial divisor of $g$ in 
$O_{d'}(\#\mathcal{G}\cdot(q^\frac{1}{2}\ln q+\ln p)\ln^{2+\gamma} p)$ deterministic time,
where $q$ is largest prime factor of $\#G_1$ and $\gamma$ any positive number.
\begin{proof}
The deterministic Pollard-Strassen \cite{pollard} algorithm can be used to find the complete 
factorization of the $q$-smooth part of $p-1$ in the stated time. The rest of the proof is based on
obvious modifications of the proofs of corollary 4.4, theorem 6.6 and on remark 4.3 from \cite{zralek}.
\end{proof}
\end{theorem}

\section{Proof of theorem \ref{main}}

If $p\le\max\left(l, [\mathcal{O}_K:\mathbb{Z}[\theta]], 
p_0, \left(1-{p_0}^{-1}\right)^\frac{-d\delta-1}{2\delta-\frac{d}{c}-\epsilon}\right)$,
where $p_0$ is the constant from lemma \ref{baza}, then we can find efficiently the complete 
factorization of $\tilde{f}$ using the deterministic Berlekamp algorithm for example. Now assume
that the reverse inequality holds. We first compute the squarefree, distinct-degree factorization of 
$\tilde{f}$ modulo $p$, that is the products $t_e$, $e\in\mathbb{N}$, of all distinct, degree $e$ irreducible divisors of $\tilde{f}$ modulo $p$.
Fix $e$; the complete factorization of $t_e$ will be found by using the following inductive procedure.
Let $g$ be a factor of $t_e$ of degree $d'$, say $d'=ke$. Suppose that $k\ge 2$.
We show below how to split $g$ nontrivially. Keep the notation of lemma \ref{baza}.
Define $\mathcal{B}_g=\pi_g\kappa\pi(\mathcal{U}\cup\mathcal{A})$.
We can assume that $\mathcal{B}_g\subset R_g^*\cup\{0\}$; in the contrary case $(b, g)$ is a
nontrivial divisor of $g$ for some $b\in\mathcal{B}_g$. Let $\mathcal{F}=\mathcal{B}_g\setminus\{0\}$. With $G$ as in theorem \ref{PH}, let $\mathcal{G}$ be the image of $\mathcal{F}$
under the homomorphism $\sigma:~R_g^*\rightarrow G$ raising every element to the power 
$\frac{p^e-1}{p-1}$. Write $G$ as $G_1\oplus G_2$ with $\#G_1=S^k$ - this condition uniquely determines $G_1$ and $G_2$. 
We apply the algorithm from theorem \ref{PH} to check whether the group $\langle\mathcal{G}\rangle_{G_1}$ is cyclic. Suppose that it is, for otherwise we would obtain a nontrivial factor of $g$. Then the order of $\langle\mathcal{G}\rangle_{G_1}$ divides $p-1$.
We will estimate this order from below to obtain a contradiction. We have 
$\#\langle\mathcal{G}\rangle_{G_1}=
\#\frac{\langle\mathcal{G}\rangle_G}{\langle\mathcal{G}\rangle_G\cap G_2}$.
The kernel of $\sigma$ has $\left(\frac{p^e-1}{p-1}\right)^k$ elements, hence
$\#\langle\mathcal{G}\rangle_G\ge
\left(\frac{p-1}{p^e-1}\right)^k\cdot\#\langle\mathcal{F}\rangle_{R_g^*}$. We appeal to lemma \ref{baza} to deduce that $\#\langle\mathcal{F}\rangle_{R_g^*}> p^{d'-\frac{d}{c}-\epsilon}-1$.
Since $S\ge (p-1)^{\frac{1}{2}+\delta}$, it follows that
$\#(\langle\mathcal{G}\rangle_G\cap G_2)\le\#G_2\le (p-1)^{\frac{k}{2}-k\delta}$.
Therefore 
$\#\langle\mathcal{G}\rangle_{G_1}>(p-1)^{\frac{k}{2}+k\delta}\cdot
\frac{p^{d'-\frac{d}{c}-\epsilon}-1}{(p^e-1)^k}$. The right hand side of this inequality is easily seen to be no less than $p-1$, which gives the desired contradiction. This means that a nontrivial divisor of $g$ 
had to be found at some stage.
Once we have completely factored $\tilde{f}$ in $\mathbb{Z}_p[Y]$, we get the complete factorization
of $f$ in $\mathbb{Z}_p[X]$ by the change of variable $Y=lX$.
\\
Obviously, the most time-consuming part of the described procedure is testing whether 
$\langle\mathcal{G}\rangle_{G_1}$ is cyclic. The stated running time follows from theorem \ref{PH},
as $\#\mathcal{G}=O_d(\ln^{ch} p)$.

\section{Concluding remarks}
\label{remarks}

Theorem \ref{PH} can be applied to get a short proof of Shoup's result from \cite{shoup2},
which states that under the ERH there is an algorithm that completely factors \textit{any} degree $d$ polynomial $f$ in $\mathbb{Z}_p[X]$ in $P^+(p-1)^\frac{1}{2}(d\ln p)^{O(1)}$ deterministic time. 
It suffices to iterate the following procedure. 
Let $g$ be a reducible factor of $f$. Adopt the notation of theorem \ref{PH}. As in von zur Gathen's algorithm, either we find directly a nontrivial divisor of $g$, or compute an element $a$ 
of $G\setminus\mathbb{Z}_p$ whose order is the power of some prime $s$ (cf. \cite{bachshallit}). 
Then we find an $s$-th power nonresidue $b$ modulo $p$. Finally we use the algorithm from theorem 
\ref{PH} with $G_1=G$ and $\mathcal{G}=\{a, b\}$ to find a nontrivial divisor of $g$. All the required 
steps can be done in the stated time.

Instead of completely factoring we could be simply interested in splitting the polynomial $f$ modulo $p$.
To this end, the running time bound obtained in theorem \ref{main} could be in some cases largely
improved. For example, if the degree of $f$ is odd then it would be sufficient to take the integer $q$
therein as the least prime such that the $q$-smooth part of $p-1$ is no less than 
$(p-1)^{\frac{1}{3}+\delta}$ for some $\delta>0$.

As another example, consider the polynomial $f=X^s-a$, where $s$ is an odd prime number and the integer $a$ is not an $s$-th power. Suppose that $f$ splits modulo $p$ into distinct linear factors, or equivalently $a$ is an $s$-th power residue modulo $p$ and $s$ divides $p-1$. In order to split $f$ modulo $p$ within the time bound of theorem \ref{main} it is then enough to choose $q$ as the least prime such that the $q$-smooth part of $p-1$ is no less than $(p-1)^{\frac{1}{s}+\delta}$ for some $\delta>0$. The point is that a nontrivial factor of $f$ modulo $p$ leads to an $s$-th root of $a$ modulo 
$p$ (cf. \cite{sze}). The remaining $s$-th roots of $a$ and hence the complete factorization of $f$ modulo $p$ can be found by computing a primitive $s$-th root of unity modulo $p$. This in turn can be done using Pila's algorithm \cite{pila}.



\bibliographystyle{amsplain}

\end{document}